\documentclass[11pt]{amsart}

\usepackage{amsmath,amssymb,amsthm,amsrefs,a4wide}
\usepackage{latexsym}
\usepackage{indentfirst}
\usepackage{graphicx}
\usepackage{placeins}
\usepackage{booktabs}
\usepackage{algorithm}
\usepackage{algorithmic}
\usepackage{multirow}

\theoremstyle{plain}

\theoremstyle{definition}

\newtheorem{example}{Example}

\theoremstyle{remark}

\newcommand{\bbm}{\begin{bmatrix}}
\newcommand{\ebm}{\end{bmatrix}}

\newcommand{\pV}{p_{\mathrm{V}}}
\newcommand{\pE}{p_{\mathrm{E}}}
\newcommand{\pVE}{p_{\mathrm{VE}}}

\newcommand{\Ber}{\mathrm{Ber}}

\begin{document}

\title[Annealed importance sampling for Ising models with mixed boundary conditions]{Annealed
  importance sampling for Ising models with mixed boundary conditions}

\author[]{Lexing Ying} \address[Lexing Ying]{Department of Mathematics, Stanford University,
  Stanford, CA 94305} \email{lexing@stanford.edu}

\thanks{This work is partially supported by NSF grant DMS 2011699.}



\begin{abstract}
  This note introduces a method for sampling Ising models with mixed boundary conditions. As an
  application of annealed importance sampling and the Swendsen-Wang algorithm, the method adopts a
  sequence of intermediate distributions that keeps the temperature fixed but turns on the boundary
  condition gradually. The numerical results show that the variance of the sample weights is
  relatively small.
\end{abstract}


\maketitle

\section{Introduction}\label{sec:intro}

This note is concerned with the Monte Carlo sampling of Ising models with mixed boundary
conditions. Consider a graph $G=(V,E)$ with the vertex set $V$ and the edge set $E$. Assume that
$V=I \cup B$, where $I$ is the subset of interior vertices and $B$ the subset of boundary
vertices. Throughout the note, we use $i,j$ to denote the vertices in $I$ and $b$ for the vertices
in $B$. In addition, $ij\in E$ denotes an edge between two interior vertices $i$ and $j$, while
$ib\in E$ denotes an edge between an interior vertex $i$ and a boundary vertex $b$. The boundary
condition is specified by $f=(f_b)_{b\in B}$ with $f_b=\pm 1$. The Ising model with the boundary
condition $f$ is the probability distribution $\pV(\cdot)$ over the configurations $s=(s_i)_{i\in
  I}$ of the interior vertex set $I$:
\begin{equation}
  \pV(s) \sim \exp\left( \beta \sum_{ij\in E} s_i s_j + \beta \sum_{ib\in E} s_i f_b \right).
  \label{eq:pV}
\end{equation}
A key feature of these Ising models is that, for certain mixed boundary conditions, below the
critical temperature the Gibbs distribution exhibits on the macroscopic scale different profiles.

\begin{figure}[h!]
  \centering
  \begin{tabular}{cc}
    \includegraphics[scale=0.55]{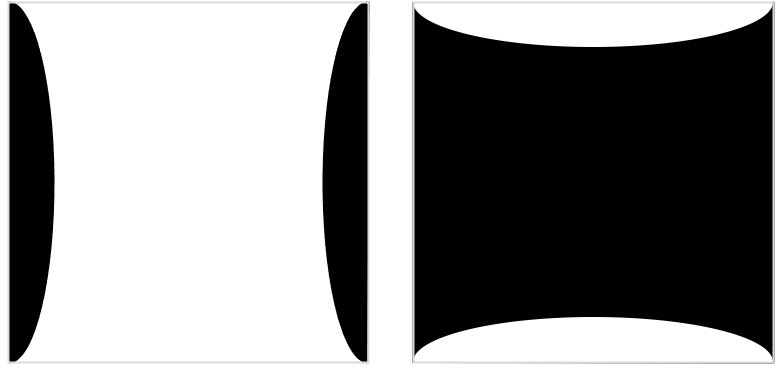} &
    \includegraphics[scale=0.55]{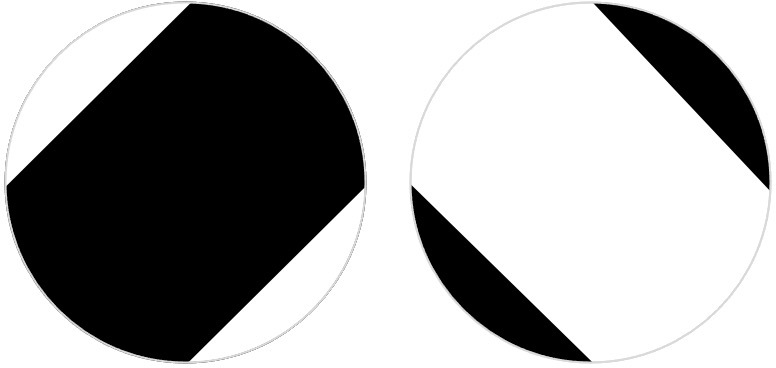} \\
    (a) & (b)
  \end{tabular}
  \caption{Ising models with mixed boundary conditions. (a) a square model and (b) a model support
    on a disk. In each case, a mixed boundary condition is specified and the model exhibits two
    dominant profiles on the macroscopic scale.}
  \label{fig:front}
\end{figure}

Figure \ref{fig:front} showcases two such examples. On the left, the square Ising lattice has the
$+1$ condition on the vertical sides but the $-1$ condition on the horizontal sides. The two
dominant profiles are a $-1$ cluster linking two horizontal sides and a $+1$ cluster linking two
vertical sides. On the right, the Ising lattice supported on a disk has the $+1$ condition on two
disjoint arcs and the $-1$ condition on the other two. The two dominant profiles are shown in Figure
\ref{fig:front}(b). Notice that in each case, the two dominant profiles have comparable
probability. Hence it is important for any sampling algorithm to visit both profiles frequently.

One of the most well-known methods for sampling Ising models is the Swendsen-Wang algorithm
\cite{swendsen1987nonuniversal}, which will be briefly reviewed in Section \ref{sec:sw}.  For Ising
models free boundary condition for example, the Swendsen-Wang algorithm exhibits rapid mixing for
all temperatures. However, for the mixed boundary conditions shown in Figure \ref{fig:front}, the
Swendsen-Wang algorithm experiences slow convergence under the critical temperatures, i.e., $T <
T_c$ or equivalently $\beta>\beta_c$ in terms of the inverse temperature. The reason is that, for
such a boundary condition, the energy barrier between the two dominant profiles is much higher than
the energies of these profiles. In other words, the Swendsen-Wang algorithm needs to break a
macroscopic number of edges between aligned adjacent spins in order to transition from one dominant
profile to the other. However, breaking so many edges simultaneously is an event with exponentially
small probability when the mixed boundary condition is specified.

Annealed importance sampling is a method by Radford Neal \cite{neal2001annealed}, designed for
sampling distributions with multiple modes. The main idea is to (1) introduce an easily-to-sample
simple distribution, (2) design a sequence of temperature dependent intermediate distributions, and
(3) generate sample paths that connects the simple initial distribution and the hard target
distribution. Annealed importance sampling has been widely applied in Bayesian statistics and data
assimilation for sampling and estimating partition functions.


In this note, we address this problem by combining the Swendsen-Wang algorithm with annealed
importance sampling. The main novelty of our approach is that, instead of adjusting the temperature,
we freeze the temperature and adjust the mixed boundary condition.

{\bf Related works.} In \cite{alexander2001spectral,alexander2001spectralB}, Alexander and Yoshida
studied the spectral gap of the 2D Ising models with mixed boundary conditions. In
\cite{ying2022double}, the {\em double flip} move is introduced for models with mixed boundary
conditions that enjoy exact or approximate symmetry. When combined with the Swendsen-Wang algorithm,
it can accelerate the mixing of these Ising model under the critical temperature
significantly. However, it only applies to problem with exact or approximate symmetries, but not
more general settings.

Recently in \cite{gheissari2018effect}, Gheissari and Lubetzky studied the effect of the boundary
condition for the 2D Potts models at the critical temperature. In \cite{chatterjee2020speeding},
Chatterjee and Diaconis showed that most of the deterministic jumps can accelerate the Markov chain
mixing when the equilibrium distribution is uniform.

{\bf Contents.} The rest of the note is organized as follows. Sections \ref{sec:sw} and
\ref{sec:ais} review the Swendsen-Wang algorithm and annealed importance sampling,
respectively. Section \ref{sec:alg} describes the algorithm and provides some numerical examples.
Section \ref{sec:disc} discusses some future directions.

\section{Swendsen-Wang algorithm}\label{sec:sw}

In this section, we briefly review the Swendsen-Wang algorithm. First, notice that
\[
\pV(s) \sim
\exp\left( \beta \sum_{ij\in E} s_i s_j + \beta \sum_{ib\in E} f_b s_i \right) = 
\exp\left( \beta \sum_{ij\in E} s_i s_j + \beta \sum_{i\in I} \left(\sum_{ib\in E} f_b\right) s_i \right).
\]
Therefore, if we interpret $h_i = \sum_{ib\in E} f_b$ as an external field, one can view the mixed
boundary condition problem as a special case of the model with external field $h=(h_i)_{i\in I}$
\begin{equation}
  \pV(s) \sim \exp\left( \beta \sum_{ij\in E} s_i s_j + \beta \sum_{i\in I} h_i s_i \right).
  \label{eq:pV2}
\end{equation}
This viewpoint simplifies the presentation of the algorithm and the Swendsen-Wang algorithm is
summarized below under this setting.

The Swendsen-Wang algorithm is a Markov Chain Monte Carlo method for sampling $\pV(s)$. In each
iteration, it generates a new configuration based on the current configuration $s$ with two
substeps:
\begin{enumerate}
\item Generate an edge configuration $w=(w_{ij})_{ij\in E}$. If the spin values $s_i$ and $s_j$ are
  different, set $w_{ij}=0$. If $s_i$ and $s_j$ are the same, $w_{ij}$ is sampled from the Bernoulli
  distribution $\Ber(1-e^{-2\beta})$, i.e., equal to 1 with probability $1-e^{-2\beta}$ and 0 with
  probability $e^{-2\beta}$.
\item Regards all edges $ij\in E$ with $w_{ij}=1$ as linked. Compute the connected components. For
  each connected component $\gamma$, define $h_\gamma = \sum_{i\in \gamma} h_i$.  Set the spins
  $(t_i)_{i\in\gamma}$ of the new configuration $t$ to $1$ with probability $e^{\beta
    h_\gamma}/(e^{\beta h_\gamma} + e^{-\beta h_\gamma})$ and to $0$ with probability $e^{-\beta
    h_\gamma}/(e^{\beta h_\gamma} + e^{-\beta h_\gamma})$.
\end{enumerate}

Associated with \eqref{eq:pV2}, two other probability distributions are important for analyzing the
Swendsen-Wang algorithm \cite{edwards1988generalization}. The first one is the joint vertex-edge
distribution
\begin{equation}
  \pVE(s,w)\sim 
  \prod_{ij\in E} \left( (1-e^{-2\beta})\delta_{s_i=s_j}\delta_{w_{ij}=1} + e^{-2\beta}\delta_{w_{ij}=0}\right)
  \cdot e^{\beta \sum_{i\in I} h_i s_i}.
  \label{eq:pVE}
\end{equation}
The second one is the edge distribution
\begin{equation}
  \pE(w) \sim
  \prod_{w_{ij}=1}(1-e^{-2\beta}) \prod_{w_{ij}=0} e^{-2\beta} \cdot 
  \prod_{\gamma\in\mathcal{C}_w} (e^{-\beta h_\gamma}+e^{\beta h_\gamma}),
  \label{eq:pE}
\end{equation}
where $\mathcal{C}_w$ is the set of the connected components induced by $w$.

Summing $\pVE(s,w)$ over $s$ or $w$ gives the following two identities
\begin{equation}
  \sum_s \pVE(s,w) = \pE(w),\quad
  \sum_w \pVE(s,w) = \pV(s).
  \label{eq:rel}
\end{equation}
(see for example \cite{edwards1988generalization}). A direct consequence of \eqref{eq:rel} is that
the Swendsen-Wang algorithm can be viewed as a data augmentation method \cite{liu2001monte}: the
first substep samples the edge configuration $w$ conditioned on the spin configuration $s$, while
the second substep samples a new spin configuration conditioned on the edge configuration $w$.

This equalities in \eqref{eq:rel} also imply that Swendsen-Wang algorithm satisfies the {\em
  detailed balance}. To see this, let us fix two spin configurations $s$ and $t$ and consider the
the transition between them. Since such a transition in the Swendsen-Wang move happens via an edge
configuration $w$, it is sufficient to show
\[
\pV(s) P_w(s,t) = \pV(t) P_w(t,s)
\]
for any compatible edge configuration $w$, where $P_w(s,t)$ is the transition probability from $s$
to $t$ via $w$. Since the transition probabilities from $w$ to the spin configurations $s$ and $t$
are proportional to $e^{\beta \sum_i h_i s_i}$ and $e^{\beta \sum_i h_i t_i}$ respectively, it
reduces to showing
\begin{equation}
  \pV(s) P(s,w) e^{\beta \sum_i h_i t_i}= \pV(t) P(t,w)e^{\beta \sum_i h_i s_i},
  \label{eq:proof}
\end{equation}
where $P(s,w)$ is the probability of obtaining the edge configuration $w$ from $s$. Using
\eqref{eq:pV2}, this is equivalent to
\begin{equation}
  e^{\beta \sum_{ij\in E} s_i s_j + \beta \sum_{i\in I} h_i s_i} P(s,w) e^{\beta \sum_{i\in I} h_i t_i} =
  e^{\beta \sum_{ij\in E} t_i t_j + \beta \sum_{i\in I} h_i t_i} P(t,w) e^{\beta \sum_{i\in I} h_i s_i}.
  \label{eq:pf}
\end{equation}
The next observation is that
\begin{equation}
  e^{\beta \sum_{ij\in E} s_i s_j} P(s,w) = e^{\beta \sum_{ij\in E} t_i t_j} P(t,w),
  \label{eq:indep}
\end{equation}
i.e., independent of the spin configuration, as explained below. First, if an edge $ij\in E$ has
configuration $w_{ij}=1$, then $s_i=s_j$. Second, if $ij\in E$ has configuration $w_{ij}=0$, then
$s_i$ and $s_j$ can either be the same or different. In the former case, the contribution to
$e^{\beta \sum_{ij\in E} s_i s_j} P(s,w)$ from $ij$ is $e^{2\beta} \cdot e^{-2\beta}=1$ up to a
uniform normalization constant. In the latter case, the contribution is also $1\cdot 1=1$ up to the
same uniform constant. After canceling the two terms of \eqref{eq:indep} in \eqref{eq:pf}, proving
\eqref{eq:proof} is equivalent to $e^{\beta \sum_{i\in I} h_i s_i} \cdot e^{\beta\sum_{i\in I} h_i t_i} =
e^{\beta \sum_{i\in I} h_i t_i} \cdot e^{\beta \sum_{i\in I} h_i s_i}$, which is trivial.

The Swendsen-Wang algorithm unfortunately does not encourage transitions between the dominant
profiles shown for example in Figure \ref{fig:front}. With these mixed boundary condition, such a
transition requires breaking a macroscopic number of edges between aligned adjacent spins, which has
an exponentially small probability.

\section{Annealed importance sampling}\label{sec:ais}

Given a target distribution $p(s)$ that is hard to sample directly, annealed importance sampling
(AIS), proposed by Neal in \cite{neal2001annealed}, introduces a sequence of distributions
\[
p_0(\cdot), \ldots, p_L(\cdot)\equiv p(\cdot),
\]
where $p_0(\cdot)$ is easy to sample and each $p_l(\cdot)$ is associated with a detailed balance
Markov Chain $T_l(s,t)$, i.e.,
\begin{equation}
  p_l(s) T_l(s,t) = p_l(t) T_l(t,s).
  \label{eq:db}
\end{equation}
The detailed balance condition can be relaxed, though it simplifies the description. Given this
sequence of intermediate distributions, AIS proceeds as follows.
\begin{enumerate}
\item Sample a configuration $s_{1/2}$ from $p_0(\cdot)$.
\item For $l=1,\ldots,L-1$, take one step (or a few steps) of $T_i(\cdot,\cdot)$ (associated with
  the distributions $p_l(\cdot)$) from $s_{l-1/2}$. Let $s_{l+1/2}$ be the resulting configuration.
\item Set $s := s_{L-1/2}$.
\item Set the weight
  \[
  w := \frac{p_1(s_{1/2})}{p_0(s_{1/2})} \cdot \cdots \cdot  \frac{p_L(s_{L-1/2})}{p_{L-1}(s_{L-1/2})}.
  \]
\end{enumerate}
The claim is that the configuration $s$ with weight $w$ samples the target distribution
$p_L(\cdot)$. To see this, consider the path $(s_{1/2}, \ldots, s_{L-1/2})$. This path is generated
with probability
\[
p_0(s_{1/2}) T_1(s_{1/2},s_{3/2}) \ldots T_{L-1}(s_{L-3/2},s_{L-1/2}).
\]
Multiplying this with $w$ and using the detailed balance \eqref{eq:db} of $T_l$ gives
\begin{align*}
  & p_0(s_{1/2}) T_1(s_{1/2},s_{3/2}) \ldots T_{L-1}(s_{L-3/2},s_{L-1/2}) \cdot \frac{p_1(s_{1/2})}{p_0(s_{1/2})} \cdot \cdots \cdot  \frac{p_L(s_{L-1/2})}{p_{L-1}(s_{L-1/2})}\\
  =&p_L(s_{L-1/2}) T_{L-1}(s_{L-1/2},s_{L-3/2}) \ldots T_1(s_{3/2},s_{1/2}),
\end{align*}
which is the probability of going backward, i.e., starting from a sample $s_{L-1/2}$ of
$p_L(\cdot)=p(\cdot)$. Taking the margin of the last slot $s_{L-1/2}$ proves that $s:=s_{L-1/2}$
with weight $w$ samples with distribution $p_L(\cdot)=p(\cdot)$.

\section{Algorithm and results}\label{sec:alg}

Our proposal is to combine AIS with the Swendsen-Wang algorithm for sampling Ising models with mixed
boundary conditions. The key ingredients are:
\begin{itemize}
\item Set the initial $p_0(\cdot)$ to be
  \[
  p_0(s) \sim \exp\left( \beta \sum_{ij\in E} s_i s_j\right).
  \]
  This initial distribution has no external field and hence can be sampled efficiently with the
  Swendsen-Wang algorithm.
\item Choose a monotone sequence of $(\theta_l)_{0\le l \le L}$ with $\theta_0=0$ and $\theta_L=1$ and set 
  at level $l$
  \[
  p_l(s) \sim \exp\left( \beta \sum_{ij\in E} s_i s_j + \beta \sum_{i\in I} (\theta_l h_i) s_i \right),
  \]
  the distribution with external field $\theta_h h$. The Markov transition matrix $T_l(\cdot,\cdot)$
  is implemented with the Swendsen-Wang algorithm associated with $p_l(\cdot)$. As proven in Section
  \ref{sec:sw}, $T_l(\cdot,\cdot)$ satisfies the detailed balance,
\end{itemize}

Below we demonstrate the performance of the proposed method with several examples. In each example,
$K=500$ samples $(s^{(k)},w^{(k)})_{1\le k\le K}$ are generated. For each sample $s^{(k)}$, the
initial choice $s^{(k)}_{1/2}$ is obtained by running $100$ iterations of the Swendsen-Wang
algorithm at $p_0(\cdot)$. In our implementation, the monotone sequence $(\theta_l)_{0\le l \le L}$
is chosen to be an equally spaced sequence with $L=400$. Although the equally-spaced sequence is not
necessarily the ideal choice in terms of variance minimization, it seems to work reasonably well for
the examples studied here.

In order to monitor the variance of the algorithm, we record the weight history at each level $l$:
\[
w^{(k)}_l = \frac{p_1(s^{(k)}_{1/2})}{p_0(s^{(k)}_{1/2})} \cdot \cdots \cdot \frac{p_l(s^{(k)}_{l-1/2})}{p_{l-1}(s^{(k)}_{l-1/2})}
\]
for $l=1,\ldots,L$. These weights are then normalized at each level $l$
\[
\tilde{w}^{(k)}_l = \frac{1}{K} \sum_{k=1}^K w^{(k)}_l.
\]
Following \cite{neal2001annealed}, we report the variance of the logarithm of the normalized weights
$\text{Var}[\{\log \tilde{w}^{(k)}_l\}]$ as a function of level $l=1,\ldots,L$. The sample
efficiency, a quantity between $0$ and $1$, is measured as
$(1+\text{Var}[\{\tilde{w}^{(k)}_L\}])^{-1}$ at level $L$.

\begin{figure}[h!]
  \centering
  \begin{tabular}{cc}
    \includegraphics[scale=0.3]{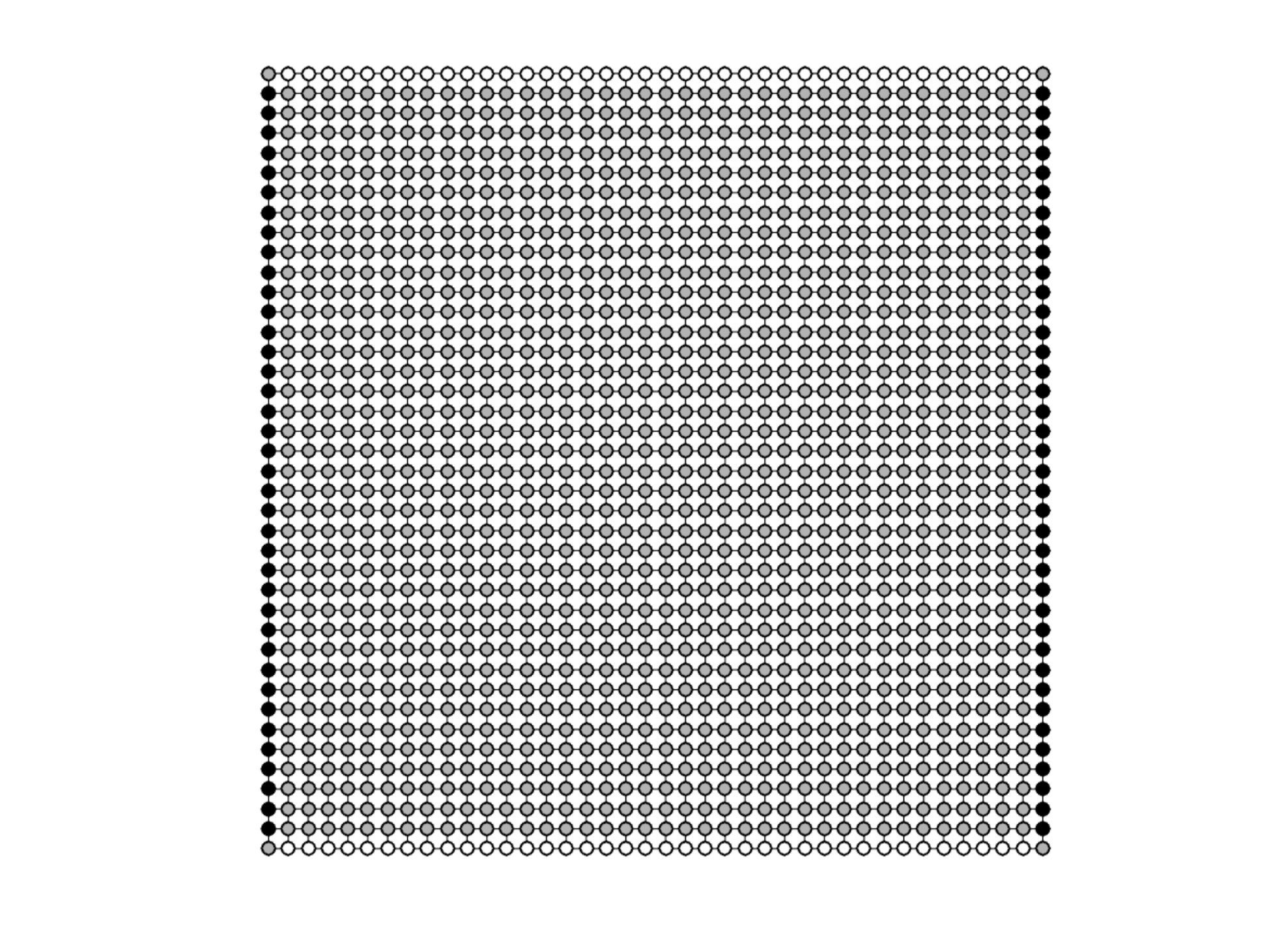} &
    \includegraphics[scale=0.3]{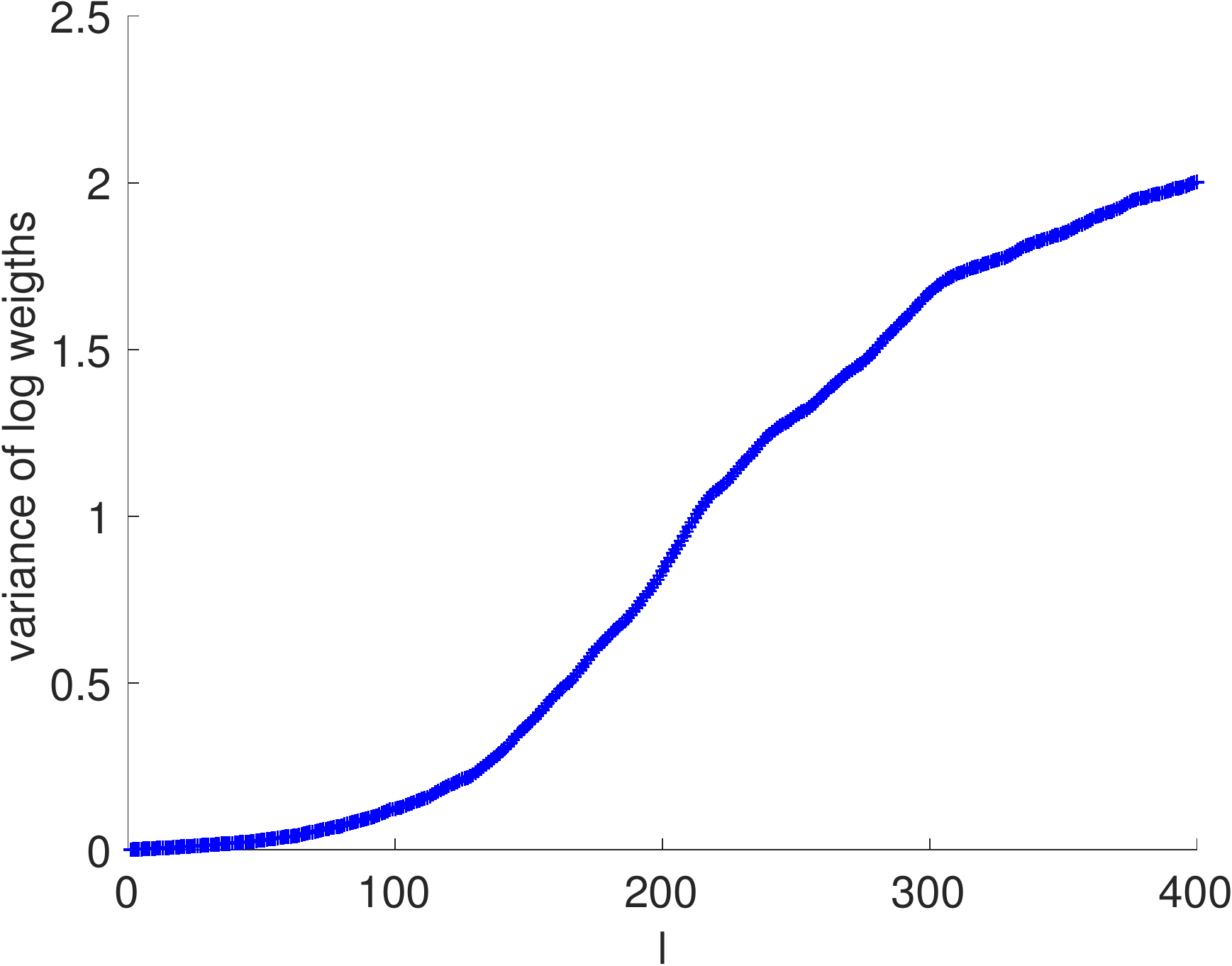} \\
    (a) & (b) 
  \end{tabular}
  \caption{(a) The lattice along with the external field. (b) The variance of the logarithm of the
    normalized weights as a function of level $l$.}
  \label{fig:ex1}
\end{figure}

\begin{figure}[h!]
  \centering
  \begin{tabular}{cc}
    \includegraphics[scale=0.3]{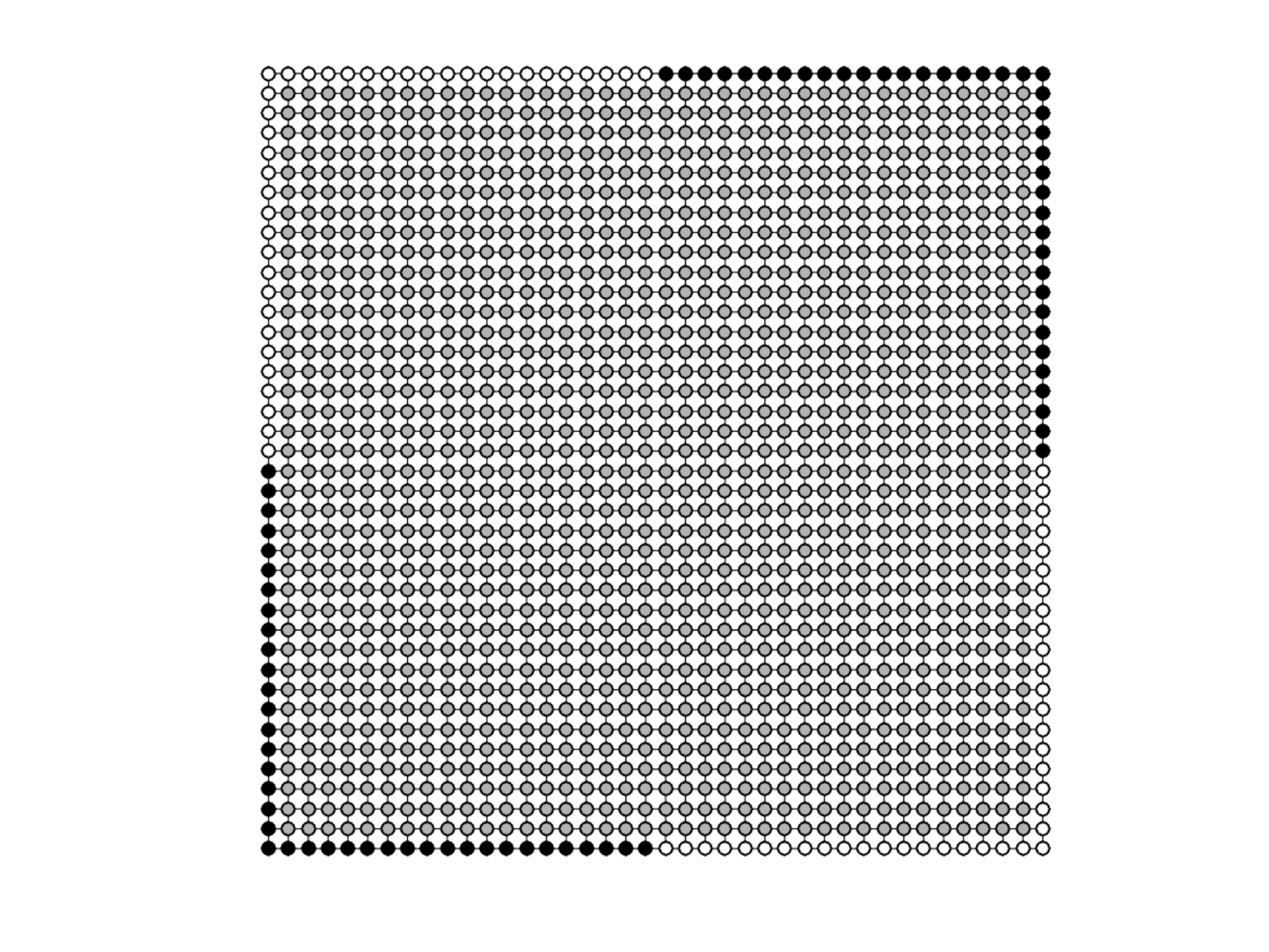} &
    \includegraphics[scale=0.3]{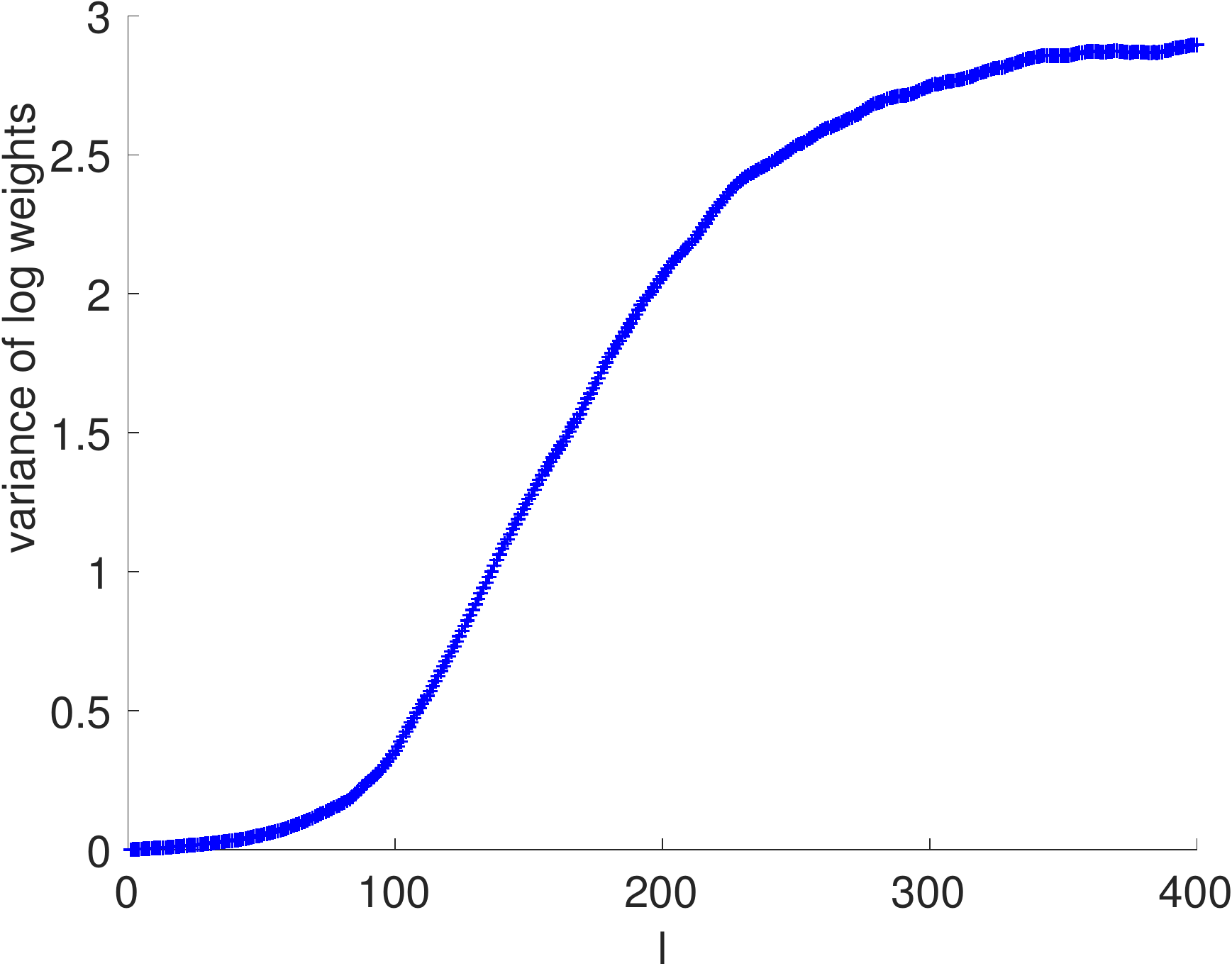} \\
    (a) & (b) 
  \end{tabular}
  \caption{(a) The lattice along with the external field. (b) The variance of the logarithm of the
    normalized weights as a function of level $l$.}
  \label{fig:ex2}
\end{figure}

\begin{example}
  The Ising model is a square lattice, as shown in Figure \ref{fig:ex1}(a). The mixed boundary
  condition is the $+1$ at the two vertical sides and $-1$ at the two horizontal sides. The
  experiments are performed for the problem size $n_1=n_2=40$ at the inverse temperature
  $\beta=0.5$. Figure \ref{fig:ex2}(b) plots the variance of the logarithm of the normalized
  weights, $\text{Var}[\{\log \tilde{w}^{(k)}_l\}]$, as a function of the level $l$. The sample
  efficiency $(1+\text{Var}[\{\tilde{w}^{(k)}_L\}])^{-1}$ is $0.26$, which translates to
  $L\cdot(1+\text{Var}[\{\tilde{w}^{(k)}_L\}])^{-1} \approx 1530$ Swendsen-Wang iterations per
  effective sample.
\end{example}

\begin{example}
  The Ising lattice is again a square as shown in Figure \ref{fig:ex2}(a). The mixed boundary
  condition is $+1$ in the first and third quadrants but $-1$ in the second and fourth
  quadrants. The experiments are performed for the problem size $n_1=n_2=40$ at the inverse
  temperature $\beta=0.5$. Figure \ref{fig:ex2}(b) plots $\text{Var}[\{\log \tilde{w}^{(k)}_l\}]$ as
  a function of the level $l$, which remain quite small. The sample efficiency
  $(1+\text{Var}[\{\tilde{w}^{(k)}_L\}])^{-1}$ is $0.09$, which translates to about
  $L\cdot(1+\text{Var}[\{\tilde{w}^{(k)}_L\}])^{-1}\approx 4470$ Swendsen-Wang iterations per
  effective sample.
\end{example}

\begin{figure}[h!]
  \centering
  \begin{tabular}{cc}
    \includegraphics[scale=0.3]{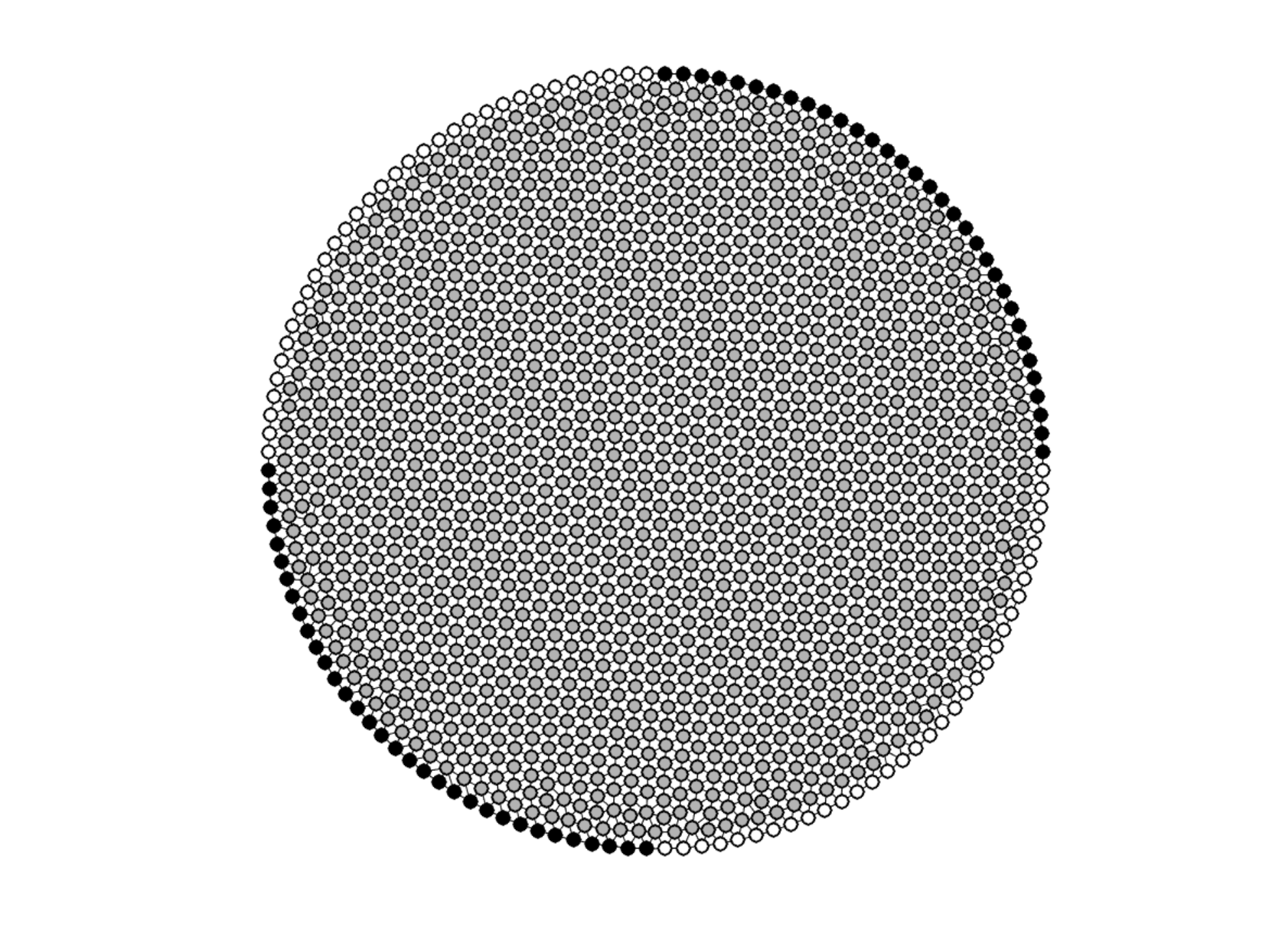} &
    \includegraphics[scale=0.3]{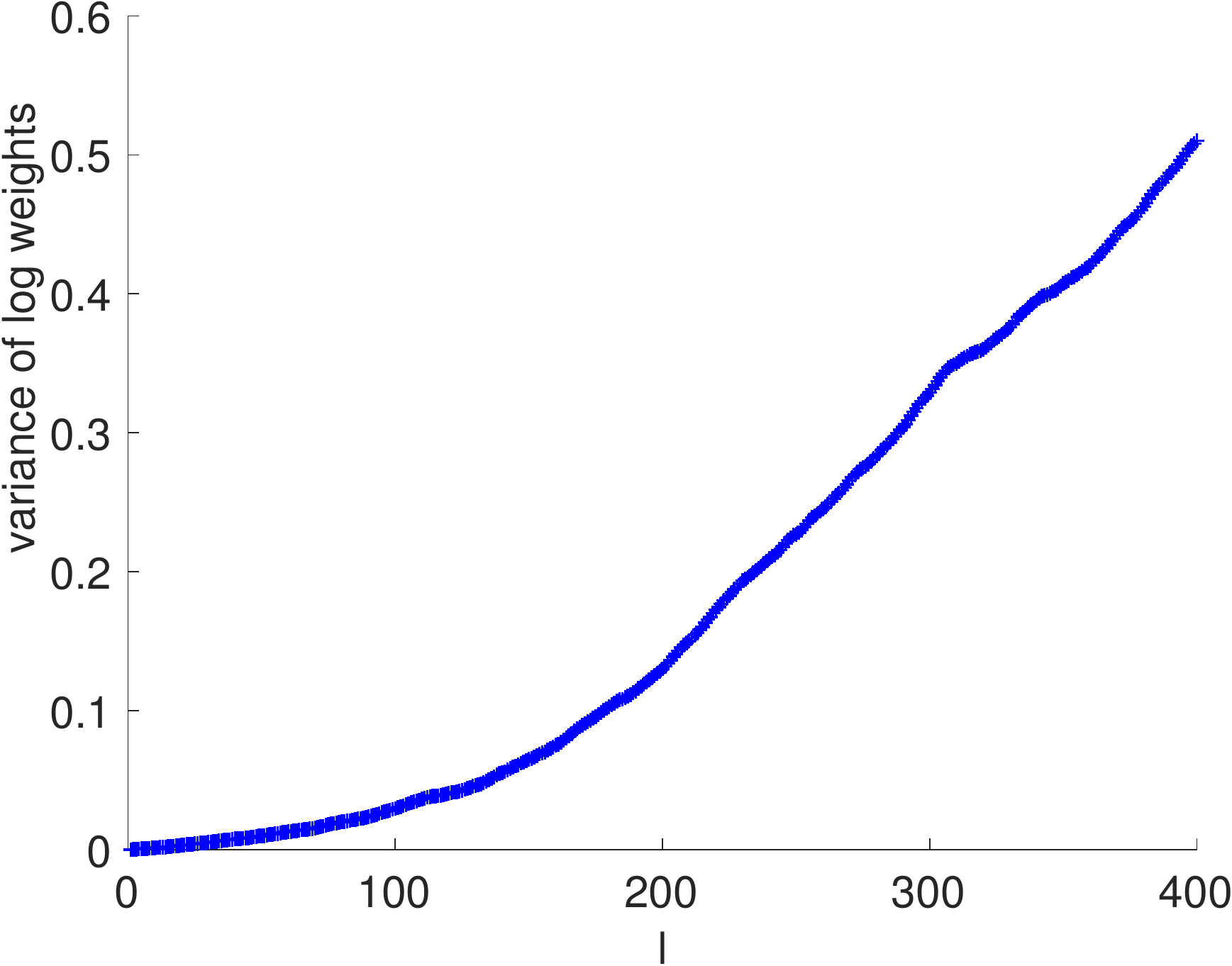} \\
    (a) & (b) 
  \end{tabular}
  \caption{(a) The lattice along with the external field. (b) The variance of the logarithm of the
    normalized weights as a function of level $l$.}
  \label{fig:ex3}
\end{figure}

\begin{example}
  The Ising model is a random quasi-uniform triangular lattice supported on the unit disk, as shown
  in Figure \ref{fig:ex3}(a). The lattice does not have rotation and reflection symmetry due to the
  random triangulation. The mixed boundary condition is $+1$ in the first and third quadrants but
  $-1$ in the second and fourth quadrants.  The experiments are performed with a finer triangulation
  with mesh size $h=0.05$ at the inverse temperature $\beta=0.3$. Figure \ref{fig:ex3}(b) plots
  $\text{Var}[\{\log \tilde{w}^{(k)}_l\}]$ as a function of the level $l$, which remain quite
  small. The sample efficiency $(1+\text{Var}[\{\tilde{w}^{(k)}_L\}])^{-1}$ is $0.67$, which
  translates to about $L\cdot(1+\text{Var}[\{\tilde{w}^{(k)}_L\}])^{-1}\approx 600$ Swendsen-Wang
  iterations per effective sample.
\end{example}

\begin{figure}[h!]
  \centering
  \begin{tabular}{cc}
    \includegraphics[scale=0.3]{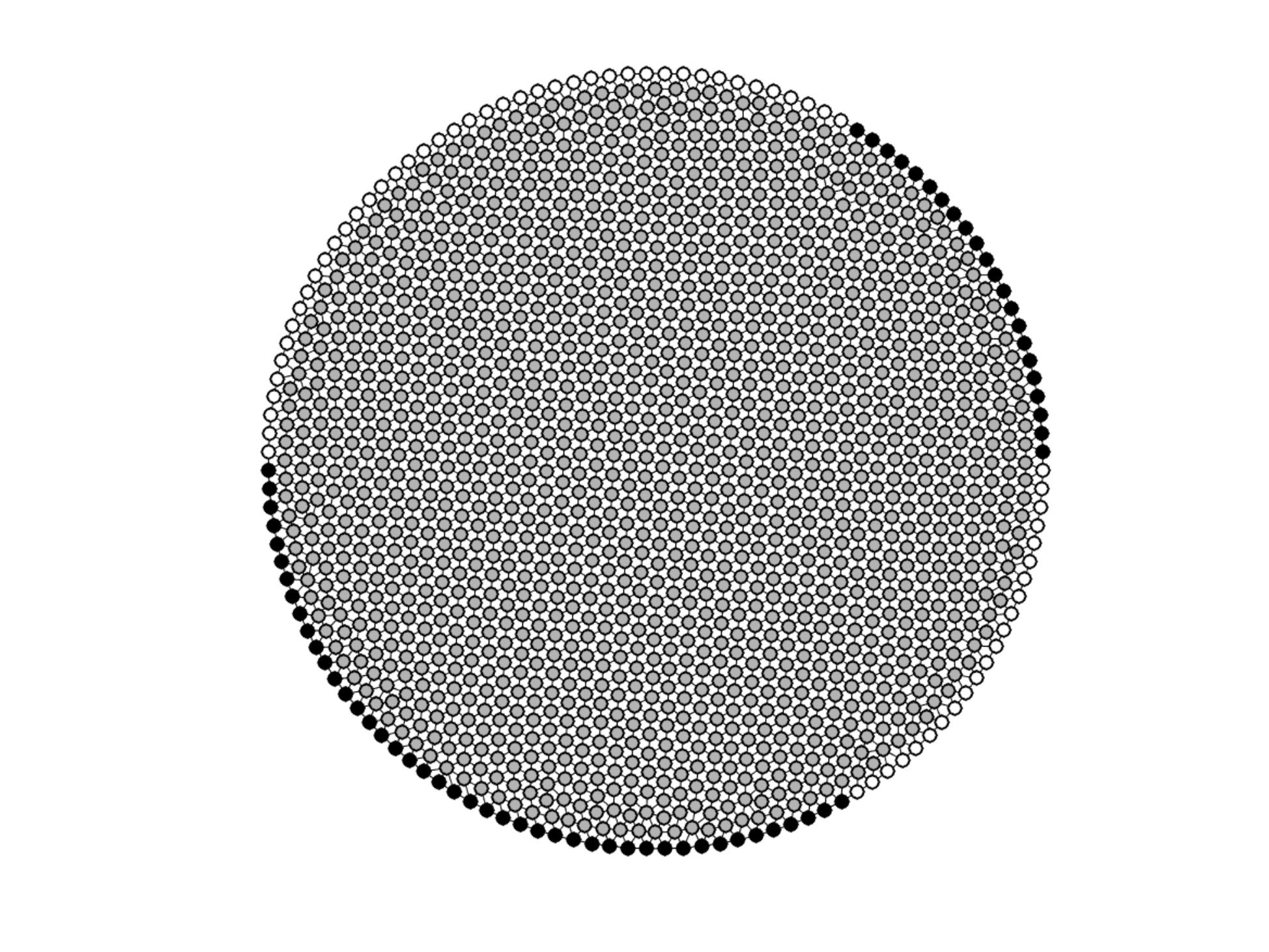} &
    \includegraphics[scale=0.3]{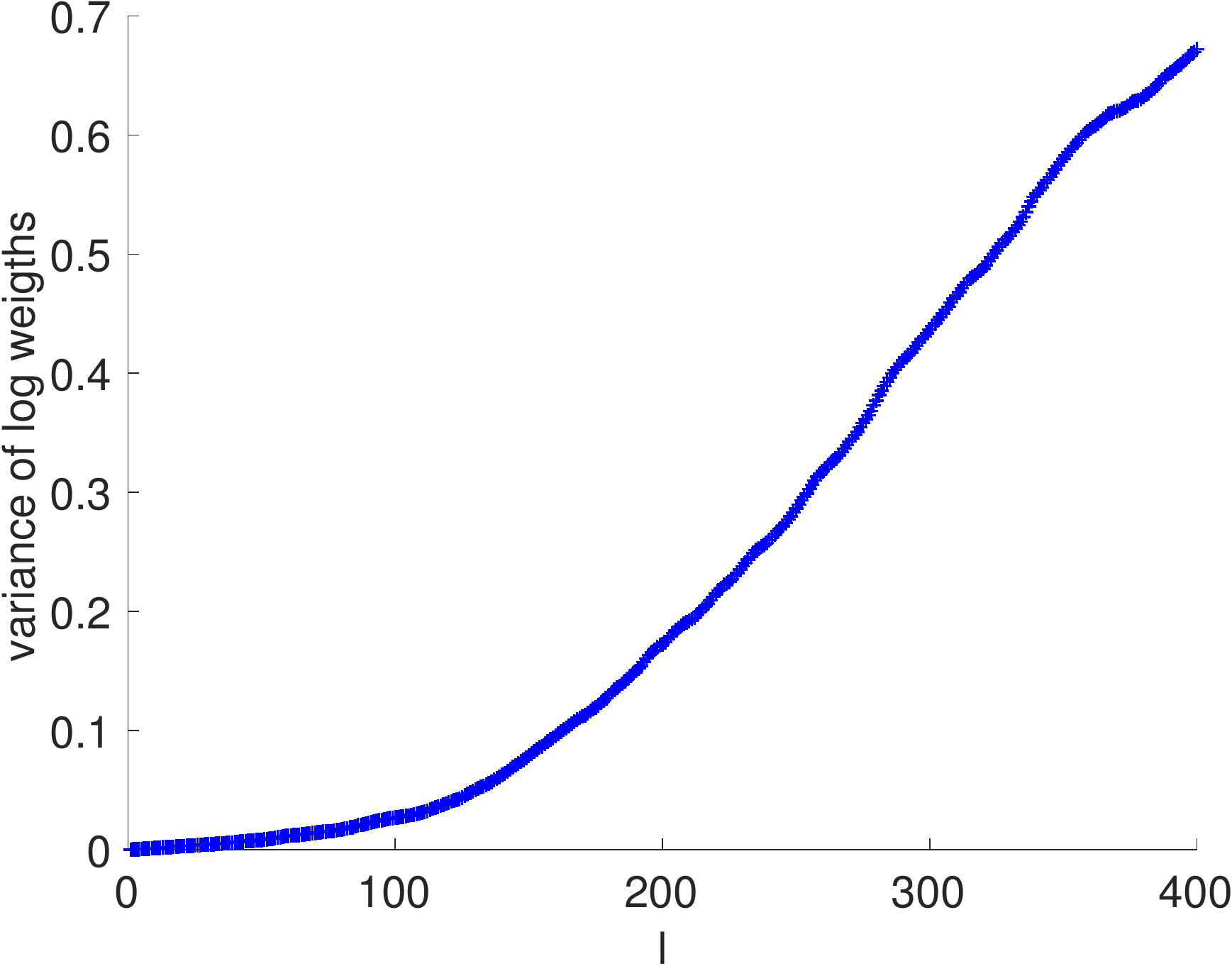} \\
    (a) & (b) 
  \end{tabular}
  \caption{(a) The lattice along with the external field. (b) The variance of the logarithm of the
    normalized weights as a function of level $l$.}
  \label{fig:ex4}
\end{figure}

\begin{example}
  The Ising model is again a random quasi-uniform triangular lattice supported on the unit disk, as
  shown in Figure \ref{fig:ex4}(a).  The mixed boundary condition is $+1$ on the two arcs with angle
  in $[0,\pi/3]$ and $[\pi,5\pi/3]$ but $-1$ on the remaining two arcs.  The experiments are
  performed with a finer triangulation with mesh size $h=0.05$ at the inverse temperature
  $\beta=0.3$. Figure \ref{fig:ex4}(b) plots $\text{Var}[\{\log \tilde{w}^{(k)}_l\}]$ as a function
  of the level $l$, which remain quite small. The sample efficiency
  $(1+\text{Var}[\{\tilde{w}^{(k)}_L\}])^{-1}$ is $0.55$, which translates to about
  $L\cdot(1+\text{Var}[\{\tilde{w}^{(k)}_L\}])^{-1}\approx 730$ Swendsen-Wang iterations per
  effective sample.
\end{example}

\section{Discussions}\label{sec:disc}

This note introduces a method for sampling Ising models with mixed boundary conditions. As an
application of annealed importance sampling with the Swendsen-Wang algorithm, the method adopts a
sequence of intermediate distributions that fixes the temperature but turns on the boundary
condition gradually. The numerical results show that the variance of the sample weights remain to be
relatively small.


There are many unanswered questions. First, the sequence of $(\theta_l)_{0\le l \le L}$ that
controls the intermediate distributions is empirically specified to be equally-spaced. Two immediate
questions are (1) what the optimal $(\theta_l)_{0\le l \le L}$ sequence is and (2) whether there is
an efficient algorithm for computing it.

Second, historically annealed importance sampling is introduced following the work of tempered
transition \cite{neal1996sampling}. We have implemented the current idea within the framework of
tempered distribution. However, the preliminary results show that it is less effective compared to
annealed importance sampling. A more thorough study is needed in this direction. 

Finally, annealed importance sampling (AIS) is a rather general framework. For a specific
application, the key to efficiency is the choice of the distribution $p_0(\cdot)$: it should be
easy-to-sample, while at the same time sufficiently close to the target distribution
$p(\cdot)$. However, since the target distribution is hard-to-sample, these two objectives often
compete with each other. There are many other hard-to-sample models in statistical mechanics. A
potential direction of research is to apply AIS with appropriate initial $p_0(\cdot)$ to these
models.

\bibliographystyle{abbrv}

\bibliography{ref}

\end{document}